\magnification=\magstep1
\input amstex
\documentstyle{amsppt}
\nologo

\pagewidth{6.00truein}
\pageheight{8.5truein}
\NoBlackBoxes
\NoRunningHeads
\topmatter
\title
Quaternionic quasideterminants and determinants
\endtitle
\author
Israel Gelfand, Vladimir Retakh, Robert Lee Wilson
\endauthor
\address
\newline
Department of Mathematics, Rutgers University,
Piscataway, NJ 08854
\endaddress
\email
\newline
I.~G.: igelfand\@math.rutgers.edu
\newline
V.~R.: vretakh\@math.rutgers.edu
\newline
R.~W.: rwilson\@math.rutgers.edu
\endemail
\endtopmatter
\document

\head Introduction \endhead

Quasideterminants of noncommutative matrices introduced in \cite
{GR, GR1} have proved to be a powerfull tool in basic problems of
noncommutative algebra and geometry (see \cite {GR, GR1-GR4,
GKLLRT, GV, EGR, EGR1, ER,KL, KLT, LST, Mo, Mo1, P, RS, RRV, Rsh,
Sch}). In general, the quasideterminants of matrix $A=(a_{ij})$
are rational functions in $(a_{ij})$'s. The minimal number of
successive inversions required to express an rational function is
called the height of this function. The "inversion height" is an
important invariant showing a degree of "noncommutativity".

In general, the height of the quasideterminants of matrices of
order $n$ equals $n-1$ (see \cite {Re}). Quasideterminants are
most useful when, as in commutative case, their height is less or
equal to $1$. Such examples include quantum matrices and their
generalizations, matrices of differential operators, etc. (see
\cite {GR, GR1-GR4, ER}).

Quasideterminants of quaternionic matrices $A=(a_{ij})$ provide a
closely related example. They can be written as a ratio of
polynomials in $a_{ij}$'s and their conjugates. In fact, our
Theorem 3.3 shows that any quasideterminant of $A$ is a sum of
monomials in the $a_{ij}$'s and the  $\bar a_{ij}$'s with real
coefficients.

The theory of quasideterminants leads to a natural definition of
determinants of square matrices over noncommutative rings.

There is a long history of attempts to develop such a theory.
These works have resulted in a number of useful generalizations of
determinants to special classes of noncommutative rings, e.g.
\cite {B, C, Ca, Di, Dy1, F, H, KS, O, R, SK, W} a.o. to
noncommutative determinants. In all these cases noncommutative
determinants can be expressed as products of quasideterminants.
It shows that quasideterminants tend to be a more fundamental
notion in noncommutative algebra.

In particular, studies of quaternionic linear algebra initiated by
Hamilton has resulted in notions of determinants due to Study and
Moore and their applications \cite {Ad, Al, As, Dy, Dy1}. We show
here that these definitions may be obtained from a general
procedure as products of quasideterminants.

A good expostion of quaternionic determinants can be found in \cite {As}.

The results of this paper can be easily generalized to a large
variety of algebras with anti-involutions.

We would like to thank S. Alesker, A. Barvinok, and I. Frenkel. Their questions
about relations between quasideterminants and determinants were very stimulating.

The second author wsa partially supported by the National Science Foundation.

\head I. Preliminary \endhead

\subhead 1.1. Quaternions \endsubhead

The algebra of quaternions $\Bbb H$ is a unital $\Bbb R$-algebra
with generators $\bold i$, $\bold j$, $\bold k $ such that
$$
\bold i^2=\bold j^2=\bold k^2=-1, \ \
\bold i \bold j=\bold k = -\bold j \bold i, \ \
\bold j \bold k=\bold i=-\bold k \bold j, \ \
\bold k \bold i=\bold j=-\bold i \bold k.
$$

Elements  $h\in \Bbb H$ can be written as $h=a+b\bold i + c\bold
j+d\bold k $ where $a,b,c,d\in \Bbb R$. Conjugation $h\mapsto
\bar h$ in $\Bbb H$ is defined by $\bar h=a-b\bold i - c\bold
j-d\bold k $ and satisfies $\overline {xy}=\bar y\bar x$. It
follows that $h\bar h=\bar h h$.

Define the norm $\nu (h)$ of a quaternion $h$ by setting $\nu
(h)=h\bar h$. The norm $\nu (h)$ is a non-negative real number
and $\nu (h)=0$ if and only if $h=0$. By definition, $|a|=\sqrt
{\nu (a)}$. It follows that $h^{-1}$ is equal to $\bar h$ up to a
central coeffcient:
$$
h^{-1}=\frac {1}{\nu (a)}\bar h .
$$

The field of real numbers $\Bbb R$ is canonically identified with
a subfield in $\Bbb H$ via embedding $a\mapsto a+0\bold i +0\bold
j + 0\bold k$. Then $h=\bar h$ if and only if $h\in \Bbb R$. The
center of $\Bbb H$ coincides with $\Bbb R$.

An embedding of the field of complex numbers $\Bbb C$ into $\Bbb
H$ is defined by an image of $\bold i\in \Bbb C$. We chose the
embedding given by $a+b\bold i \mapsto a+b\bold i +0\bold j +
0\bold k$, where $a,b\in \Bbb R$. The embedding $\Bbb C$ into
$\Bbb H$ agrees with the conjugation.

Under these identifications every quaternionic number $h $ can be
written as $h = \alpha  +\bold j\beta$ where $\alpha, \beta\in
\Bbb C$.

Let $M(n, F)$ be the algebra of $n\times n$-matrices over a field
$F$. Define a map $\theta :\Bbb H\rightarrow M(2, \Bbb C)$ by
setting
$$
\theta (h)=\left ( \matrix
\alpha & -\bar \beta \\
\beta & \bar \alpha
\endmatrix \right ).
$$
It is well-known that $\theta $ is a homomorphism of $\Bbb R$-algebras.
Also, $\det \theta (h)=\nu (h)$.

If matrix $A=(a_{ij})$, $i,j=1,\dots , n$ is a quaternionic
matrix, define its corresponds its Hermitian dual matrix to be
$A^*=(b_{ij})$ where $b_{ij}=\bar a_{j,i}$ for all $i,j$. $A$ is a
Hermitian matrix if $A=A^*$.

\subhead 1.2. Quasideterminants \endsubhead

\subhead 1.1.2. Definition of quasideterminants\endsubhead

The definition of quasideterminants and their main properties were
given in \cite {GR. GR1-GR4}.

Let $A=(a_{ij}), i\in I, j\in J$ be a matrix over a ring $R$.
Construct the following submatrices of A: submatrix $A^{ij}$,
$i\in I$, $j\in J$ obtained from $A$ by deleting its $i$-th row
and $j$-th column; row submatrix $r_i^j$ obtained from $i$-th row
of $A$ by deleting the element $a_{ij}$; column submatrix $c_j^i$
obtained from $j$-th column of $A$ by deleting the element
$a_{ij}$.

\proclaim{Definition 1.1} Let $A$ be a square matrix. We say that
the quasideterminant $|A|_{ij}$ is defined if the submatrix
$A^{ij}$ is invertible over the ring $R$. In this case
$$
|A|_{ij}=a_{ij}-r_i^j(A^{ij})^{-1}c_j^i. \tag 1.1
$$
\endproclaim

We set $|A|_{ij}=a_{ij}$ if $I=\{i\}$, $J=\{j\}$.

If inverses to the quasideterminants $|A^{ij}|_{pq}$, $p\neq i,
q\neq j$, are defined then

$$|A|_{ij} = a_{ij} -\sum
a_{iq}|A^{ij}|^{-1}_{pq} a_{pj}.\tag 1.2
$$
Here the sum is taken over all $p\in I\smallsetminus\{i\},q\in
J\smallsetminus\{j\}$.

If $A$ is an $n\times n$-matrix there exist up to $n^2$
quasideterminants of $A$.

If $A$ is a square matrix we call a quasideterminant of any
$k\times k$-submatrix  a $k$-quasiminor of $A$.

\example{Examples 1.2}\endexample

1) For a matrix $A=(a_{ij}), i,j =1,2$
there exist four quasideterminants if the corresponding entries are invertible
$$\matrix
|A|_{11} = a_{11} - a_{12}\cdot a^{-1}_{22}\cdot a_{21},\quad
&|A|_{12}=a_{12}-a_{11}\cdot a^{-1}_{21}\cdot a_{22},\\
|A|_{21} = a_{21} - a_{22}\cdot a^{-1}_{12}\cdot a_{11},\quad &
|A|_{22}=a_{22}-a_{21}\cdot a^{-1}_{11}\cdot a_{12}.\endmatrix
$$

If an entry, say $a_{22}$ is not invertible then the quasidetrminant $|A|_{11}$
does not exist.

2) For a matrix $A=(a_{ij}), i,j=1,2,3$ there exist up to 9 quasideterminants
but we will write here only
$$\matrix
|A|_{11}=a_{11}-a_{12}(a_{22}-a_{23}a_{33}^{-1}a_{32})^{-1}a_{21} &-a_{12}(a_{32}-a_{33}\cdot
a_{23}^{-1} a_{22})^{-1} a_{31}\\
 \qquad -a_{13}(a_{23}-a_{22}a_{32}^{-1}a_{33})^{-1}a_{21} &-a_{13}(a_{33}-a_{32}\cdot
a^{-1}_{22}a_{23})^{-1}a_{31}\endmatrix
$$
provided all inverses are defined.

\remark{Remark} If each $a_{ij}$ is an invertible morphism
$V_j\to V_i$ in an additive category then the quasideterminant
$|A|_{pq}$ is also a morphism from the object $V_q$ to the object
$V_p$.\endremark
\medskip

\example {Example 1.3}\endexample Quasideterminants are not
generalizations of determinants over a commutative ring but
generalizations of the ratio of two determinants. If in formulas
(1.1), or (1.2) the variables $a_{ij}$ commute then
$$
|A|_{pq} = (-1)^{p+q} \frac{\det A}{\det A^{pq}}.
$$

As a general rule, quasideterminants are rational functions in
their entries but not polynomials.

The following theorem was conjectured by I. Gelfand and V. Retakh, and
proved by Reutenauer \cite {Re}.
\proclaim{Theorem 1.4}  Quasideterminants of the $n \times n$-matrix
with formal entries have the inversion height $n-1$ over the free skew-field
generated by its entries.
\endproclaim

{\bf 1.3. General properties of quasideterminants}

\subhead Elementary properties of quasideterminants\endsubhead

\roster
\item"{i)}"  The quasideterminant $|A|_{pq}$ does not depend on the permutation
of rows and columns in the matrix $A$ if the $p$-$th$ row and the $q $-$th$
column are not changed;

\item"{ii)}"  {\it The multiplication of rows and columns.}  Let the
matrix $B$ be constructed from the matrix $A$ by multiplication of its
$i$-$th$ row by a scalar $\lambda$  from the left.  Then
$$
 |B|_{kj}=\cases \lambda |A|_{ij} \qquad&\text{ if } k = i\\
|A|_{kj} \qquad&\text{ if } k \neq i \text{ and } \lambda \text{ is
invertible.}\endcases
$$
Let the matrix $C$ be constructed from the matrix A by multiplication
of its $j$-$th$ column by a scalar $\mu$ from the right.  Then
$$
|C|_{i\ell}=\cases |A|_{ij} \mu \qquad&\text{ if } \ell = j\\
|A|_{i\ell} \qquad&\text{ if } \ell \neq j \text{ and } \mu \text{ is
invertible.}\endcases
$$

\item"{iii)}" {\it The addition of rows and columns.} Let the matrix
$B$ be constructed by adding to the $p$-th row of the matrix $A$
its $k$-th row multiplied by a scalar $\lambda$ from the left,
where $k\neq p$. Then
$$
|A|_{ij} = |B|_{ij}, \quad  i=1, \dots k-1, k+1,\dots n,
 j=1, \dots, n.
$$
Let the matrix $C$ be constructed by addition to the $q$-th column
of the matrix $A$ its $\ell$-th column multiplied by a scalar
$\lambda$ from the right, where $\ell \neq q$.  Then
$$
|A|_{ij}= |C|_{ij} , \, i=1,\dots, n, j=1,\dots ,\ell -1,\ell +1,\dots n.
$$
\endroster

\subhead Homological relations \endsubhead The ratio of two
quasideterminants of the same square matrix is the ratio of two
quasiminors.  This is a consequence of the following {\it
homological} relations \cite {GR, GR1}: \proclaim{Theorem 1.5}

a) Row homological relations:
$$
-|A|_{ij} \cdot |A^{i\ell}|^{-1}_{sj} = |A|_{i\ell}\cdot
|A^{ij}|^{-1}_{s\ell}\qquad \forall s\neq i
$$

b) Column homological relations:
$$
-|A^{kj}|^{-1}_{it} \cdot |A|_{ij} = |A^{ij}|^{-1}_{kt}\cdot |A|_{kj}
\qquad \forall r\neq j
$$
\endproclaim

\subhead Heredity \endsubhead  Let $A=(a_{ij}),
i,j=1,\dots, n$, be a matrix and
$$
A=\pmatrix A_{11} & \dots & A_{1s}\\
 {}&{} &{}\\
A_{s1} &\dots &A_{ss}\endpmatrix
$$
be a block decomposition of $A$.  Denote by $\tilde A= (A_{ij}),
i,j=1,\dots, s$, the matrix with $A_{ij}$'s as entries. Suppose
that
$$A_{pq} = \pmatrix
a_{k\ell} & \dots & a_{k,\ell+m}\\
{} &\dots &{}\\
a_{k+m,\ell} & \dots & a_{k+m,\ell+ m}\endpmatrix.
$$
is a square matrix and that $|\tilde A|_{pq}$ is defined.

\proclaim{Theorem 1.6}
$$
|A|_{k'\ell'}=\big\vert|\tilde A|_{pq}\big\vert_{k'\ell'} \text{ for }
k'=k,\dots, k+m,\
\ell'=\ell,\dots, \ell +m.
$$
\endproclaim

In particular,  let $A=\pmatrix A_{11}& A_{12}\\A_{21}&
A_{22}\endpmatrix$ be a block decomposition of $A=(a_{ij}),
i,j=1,\dots, n$.  Let $A_{11}$ be a $k\times k$-matrix and let
the matrix $A_{22}$ be invertible.

Then
$|A|_{ij}= |A_{11}-A_{12} A_{22}^{-1}\cdot A_{21} |_{ij}
\text{ for } i,j=1.\dots, k.$

In other words, the quasideterminant $|A|_{ij}$ of an $n\times n$
matrix can be computed in two steps: first, consider the
quasideterminant $|\tilde A|_{11} = A_{11} - A_{12} \cdot
A_{22}^{-1}\cdot A_{21}$ of a $2\times 2$-matrix $\tilde
A=\pmatrix A_{11} &A_{12}\\
A_{21} & A_{22}\endpmatrix$; and second, consider the corresponding
quasideterminant of the $k\times k$-matrix $|\tilde A|_{11}$.

\subhead Basic identities\endsubhead The following noncommutative
analogue of Sylvester's identity is closely related with the
heredity property.  Let $A=(A_{ij}), i,j=1,\dots, n$, be a matrix
over a noncommutative skew-field. Suppose that the submatrix
$A_0=(a_{ij}), i,j=1, \dots, k$, is invertible.  For
$p,q=k+1,\dots, n$ set
$$
b_{pq} = \vmatrix &A_0 & a_{1q}\\
&{}& \vdots\\
&{}&a_{kq}\\
&a_{p1}\dots a_{pk} &a_{pq}\endvmatrix_{pq} ;
$$
$$
B=(b_{pq}), p,q = k+1,\dots, n.
$$

\proclaim{Theorem 1.7} For $i,j = k+1,\dots, n$,
$$
|A|_{ij} = |B|_{ij}.
$$
\endproclaim

For example, a quasideterminant of an $n\times n$-matrix $A$ is
equal to the corresponding quasideterminant of a $2\times
2$-matrix consisting of $(n-1)$- quasiminors of the matrix $A$,
or to the corresponding quasideterminant of an
$(n-1)\times(n-1)$-matrix consisting of $2$-quasiminors of the
matrix $A$. One can use any of these procedures for a definition
of quasideterminants.

\head II. Noncommutative determinants \endhead

Let $A=(a_{ij})$, $i,j=1,\dots , n$, be a generic matrix over a
division algebra $R$, i.e. assume all square submatrices of $A$
are invertible. For any orderings $I=(i_1,\dots , i_n)$,
$J=(j_1,\dots , j_n)$ of $\{1,\dots , n\}$ define $A^{i_1\dots
i_k, j_1\dots j_k}$ to be submatrices of $A$ obtained by deleting
rows with indices $i_1,\dots , i_k$ and columns with indices
$j_1,\dots , j_k$. Then set
$$
D_{I,J}(A)=|A|_{i_1 j_1}|A^{i_1 j_1}|_{i_2 j_2}|A^{i_1i_2, j_1j_2}|_{i_3 j_3}
\dots a_{i_n j_n}.
$$

In commutative case all $D_{I,J}(A)$ are equal up to a sign to
the determinant of $A$. When $A$ is a quantum matrix all
$D_{I,J}(A)$ are equal up to a scalar to the quantum determinant
of $A$ \cite {GR, GR1, KL}. The same is true for other well-known
algebras. This  gives us a reason to call the $D_{I,J}(A)$ the
$(I,J)$-predeterminants of $A$. From a ``categorical point of
view" the expressions $D_{I,I'}$ when $I'=(i_2, i_3,\dots , i_n,
i_1)$ are particularly important. Set $D_I=D_{I,I'}$. There are
$n!$ expressions $D_I$. It is convinient to have a basic
prederminant $\Delta =D_{\{12\dots n\},\{12\dots n\}}$.

We use the homological relations for quasideterminants to compare
of different $D_{I,J}$'s. For example, if $I=(i_1,i_2,\dots ,
i_n)$, $J=(j_1,j_2,\dots , j_n)$ and $I'=(i_2,i_1,\dots , i_n)$,
$J'=(j_2,j_1,\dots , j_n)$ it is enough to compare $|A|_{i_1
j_1}|A^{i_1 j_1}|_{i_2 j_2}$ and $|A|_{i_2 j_2}|A^{i_2 j_2}|_{i_1
j_2}$. The homological relations show that these expressions are
in a sense ``conjugate":
$$
|A|_{i_2 j_2}|A^{i_2 j_2}|_{i_1 j_1}=|A^{i_1 j_1}|_{i_2 j_2}|A^{i_2 j_1}|_{i_1 j_2}^{-1}
|A|_{i_1 j_1}|A^{i_1 j_1}|_{i_2 j_2}^{-1}|A^{i_2 j_1}|_{i_1 j_2}|A^{i_2 j_2}|_{i_1 j_1}.
$$

Let $R^*$ be the monoid of invertible elements in $R$ and $\pi
:R^*\rightarrow R^*/[R^*,R^*]$ be the canonical homomorphism of
$R^*$. \proclaim{Proposition 2.1} $\pi (\Delta (A))$ is the
Dieudonne determinant of $A$ and $\pi (D_{I,J})=p(I)p(J)\Delta$
where $p(I)$ is the parity of ordering $I$.
\endproclaim

Note, that in \cite {Dr} P. Draxl introduced a Dieudonne
predeterminant, denoted $\delta \epsilon \tau $. For a generic
matrix $A$ over a division algebra there exists the Gauss
decomposition $A=UDL$ where $U, D, L$ are upper-triangular,
diagonal, and low-triangular matrices. Draxl denoted the product
of diagonal elements in $D$ from top to the bottom by $\delta
\epsilon \tau  (A)$ and called it the Dieudonne predeterminant of
$A$. For non-generic matrices Draxl used the Bruhat decomposition
instead of Gauss decomposition.

\proclaim{Proposition 2.2} $\delta \epsilon \tau  (A)=\Delta (A)$.
\endproclaim

\demo{Proof} Let $y_1, \dots , y_n$ be the diagonal elements in
$D$ from top to the bottom. As shown in \cite {GR1},
$y_k=|A^{12\dots k-1, 12\dots k-1}|_{kk}$. Then $\delta \epsilon
\tau  (A)=y_1y_2\dots y_n=\Delta (A)$.
\enddemo
\head III. Quaternionic quasideterminants \endhead

\subhead 3.1 Norms of quaternionic matrices \endsubhead

Let $M(n, \Bbb H)$ be the $\Bbb R$-algebra of quaternionic
matrices of order $n$. There exists a multiplicative functional
$\nu : M(n, \Bbb H)\rightarrow \Bbb R_{\geq 0}$ such that

i) $\nu (A)=\nu (a)$ if $A\in M(1, \Bbb H)$, $A=(a)$,

ii) $\nu (A)=0$ if and only if  the matrix $A$ is non-invertible,

iii)  If $A '$ is obtained from $A $ by adding a left-multiple
of a row to another row or a right-multiple of a column to another
column, then $\nu  (A ')=\nu  (A)$.

We will show later that for a generic quaternionic matrix
$A=(a_{ij})$, $i,j=1,\dots , n$, a norm can be defined by
induction on $n$ by formula $\nu (A)=\nu (|A|_{11})\nu (A^{11})$
and that $\nu (A)$ is a polynomial in $a_{ij}$'s, $\bar a_{ij}$'s.

\subhead 3.2. Quasideterminants of quaternionic matrices
\endsubhead

Let $A=(a_{ij})$ ,  $i,j=1,\dots , n$, be a quaternionic matrix.
Let $I$ and $J$ be two ordered sets of natural numbers $1\leq
i_1, i_2,\dots , i_k\leq n$ and $1\leq j_1,j_2,\dots , j_k\leq n$
such that the $i_p$'s are distinct and the $j_p$'s are distinct.
If $k=1$ set  $m_{I,J}(A)=a_{i_1j_1}$. If $k\geq 2$ set
$$
m_{I,J}(A)=a_{i_1j_2}\bar a_{i_2j_2}a_{i_2j_3}\bar a_{i_3j_3}
a_{i_3j_4}\dots \bar a_{i_kj_k}a_{i_kj_1}.
$$

If the matrix $A$ is Hermitian then
$$
m_{I,J}(A)=a_{i_1j_2}a_{j_2i_2}a_{i_2j_3} a_{j_3i_3}
a_{i_3j_4}\dots a_{j_ki_k}a_{i_kj_1}.
$$
\remark {Remark} The monomials  $m_{I,J}(A)$ are conviniently
described in terms of paths in a graph.
Let $\Gamma _n$ be the complete oriented 2-graph with vertices
$1,\dots , n$, the arrows from $p$ to $q$ being labelled  by
$a_{pq}$ and $\bar a_{qp}$. The monomial $m_{I,J}(A)$
corresponds then to a path from $i_1$ to $j_1$.
\endremark

To a quaternionic  matrix $A=(a_{pq})$, $p,q=1,\dots , n$, and to
a fixed row index $i$ and a column index $j$ we associate a
polynomial in $a_{pq}$'s, $\bar a_{pq}$'s which we call the
$(i,j)$-th double permanent of $A$. \proclaim {Definition 2.1}
The $(i,j)$-th double permanent of $A$ is the sum
$$
\pi _{ij}(A)=\sum  m_{I,J}(A),
$$
taken over all orderings $I=\{i_1,\dots , i_n\}$,  $J=\{j_1,\dots
, j_n\}$ of $\{1,\dots , n\}$ such that $i_1=i$ and  $j_1=j$ .
\endproclaim

\example{Example 3.2} For $n=2$
$$
\pi _{11}(A)=a_{12}\bar a_{22}a_{21}.
$$
For $n=3$
$$
\pi _{11}(A)=a_{12}\bar a_{32}a_{33}\bar a_{23}a_{21}+
a_{12}\bar a_{22}a_{23}\bar a_{33}a_{31}+
a_{13}\bar a_{33}a_{32}\bar a_{22}a_{21}+
a_{13}\bar a_{23}a_{22}\bar a_{32}a_{31}.
$$
\endexample

For any submatrix $B$ of $A$ denote by $B^c$ its complimentary
submatrix. The matrix $B^c$ is obtained from $A$ by omitting all
rows and columns containing elements from $B$. If $B$ is a
$k\times k$ matrix then $B^c$ is a $(n-k)\times (n-k)$ matrix.

One of our main results is the following theorem. \proclaim
{Theorem 3.3} If the quasideterminant $|A|_{ij}$ of a quaternionic
matrix is defined then
$$
\nu ( A^{ij}) |A|_{ij}= \sum (-1)^{K(B)-1}\nu ( B^c)\pi _{ij}(B)
\tag 3.1
$$
where the sum is taken of all square submatrices $B$ of $A$
containing $a_{ij}$ and $K(B)$ means the order of $B$.
\endproclaim

Note that the quasideterminant $|A|_{ij}$ is defined if the matrix
$A^{ij}$ is invertible. In this case $\nu (A^{ij})$ is invertible
and formula  (3.1) gives an expression for $|A|_{ij}$.

Note that both sides in formula (3.1) are polynomials in the
$a_{ij}$'s, the $\bar a_{ij}$'s with real coefficients. The
right-hand side in (2.1) is a linear combination  of monomials of
lengths $1,3,\dots , 2n-1$. Let $\mu (n)$ be the number of such
monomials for a matrix of order $n$. Then $\mu (n)=1+(n-1)^2\mu
(n-1)$.

\example{Example 3.4} For $n=2$
$$
\nu ( a_{22}) |A|_{11}=\nu ( a_{22})a_{11}-a_{12}\bar a_{22}a_{21}.
$$
For $n=3$
$$\matrix
\nu (A^{11})|A|_{11}=\nu (A^{11})a_{11}-\\
\qquad -\nu (a_{33})a_{12}\bar a_{22}a_{21}-
\nu (a_{23})a_{12}\bar a_{32}a_{31}-\nu (a_{32})a_{13}\bar a_{23}a_{21}
-\nu (a_{22})a_{13}\bar a_{33}a_{31}+\\
\qquad +a_{12}\bar a_{32}a_{33}\bar a_{23}a_{21}+
a_{12}\bar a_{22}a_{23}\bar a_{33}a_{31}+
a_{13}\bar a_{33}a_{32}\bar a_{22}a_{21}+
a_{13}\bar a_{23}a_{22}\bar a_{32}a_{31}.\endmatrix
$$
\endexample

Example 3.4 shows how to simplify the general formula for
quasideterminants of matrix of order $3$ (see example 1.2) for
quaternionic matrices.

Formula (3.1) gives a recurrence relation for norms of
quaternionic matrices. From this relation one can obtain an
expression of $\nu (A)$ as a polynomial in $a_{ij}$'s and their
conjugates (see also subsection 4.3 below).

\proclaim{Corolarry 3.5}
$$
\nu (A) \nu ( A^{ij}) = \sum (-1)^{K(B)-1}\nu ( B^c)\pi _{ij}(B))
(\sum (-1)^{K(B)-1}\nu ( B^c)\overline {\pi _{ij}(B)}).
$$
\endproclaim

\example{Example 3.6} For $n=2$
$$
\nu (A)=\nu (a_{11}) \nu (a_{22})+\nu (a_{12})\nu (a_{21})-
a_{12}\bar a_{22}a_{21}\bar a_{11}-
a_{11}\bar a_{21}a_{22}\bar a_{12}.
$$
\endexample

For a square quaternionic matrix $A=(a_{ij})$ set
$Q_{ij}(A)=|A|_{ij}\nu (A^{ij})$. Then one can rewrite formula
(3.1) as
$$
Q_{ij}(A)=\nu (A^{ij})a_{ij}-\sum _{p\neq i, q\neq
j}a_{iq}\overline {Q_{pq}(A^{ij})}a_{pj}.
$$
If $A$ is a Hermitian matrix then  $\bar Q_{pq}(A^{ij})=Q_{qp}(A^{ij})$ and
$$
Q_{ij}(A)=\nu (A^{ij})a_{ij}-\sum _{p\neq i, q\neq j
}a_{iq}Q_{qp}(A^{ij})a_{pj}.
$$

When $A$ is a complex matrix, formula (3.1) can be written as
$$
\det A \det \bar A^{ij}=(-1)^{i+j}(a_{ij} \det \bar A^{ij}-
\sum _{p\neq i, q\neq j}a_{iq}\det \bar A^{ij}\det A^{ip,jq}a_{pj}).
$$
This implies the standard row/column expansion of the determinant.

\head IV. Determinants of quaternionic matrices \endhead

\subhead 4.1. Dieudonne determinant  \endsubhead

Let $A=(a_{ij})$, $i,j=1,\dots , n$ be a generic quaternionic
matrix. Recall that for any orderings $I=(i_1,\dots , i_n)$,
$J=(j_1,\dots , j_n)$ of $\{1,\dots , n\}$ we defined expressions
$$
D_{I,J}(A)=|A|_{i_1 j_1}|A^{i_1 j_1}|_{i_2, j_2}|A^{i_1i_2,
j_1j_2}|_{i_3 j_3} \dots a_{i_n j_n}.
$$

In the quaternionic case the Dieudonne determinant $D$ maps
$$
D: M_n(\Bbb H)\rightarrow \Bbb R_{\geq 0}
$$
(see \cite {As}).

The following proposition generalizes a result from \cite {VP}.
\proclaim {Proposition 4.1} In the quaternionic case
$|D_{I,J}(A)|$, the absolute value of $D_{I,J}(A)$, is equal to
$D(A)$ for any $I, J$.
\endproclaim
A proof of  Proposition 4.1 follows from the homological relations
for quasideterminants.

\subhead 4.2. Moore determinants of Herimitian quaternionic
matrices \endsubhead

A quaternionic matrix $A=(a_{ij})$, $i,j=1,\dots , n$ is
Hermitian if $a_{ji}=\bar a_{ji}$ for all $i,j$. It follows that
all diagonal elements of $A$ are real numbers and that the
submatrices $A^{11}$,   $A^{12,12}$, $\dots $ are Hermitian.

A determinant of Hermitian matrices was introduced by E.M. Moore
in 1922 \cite {M, MB}. Here we recall the original definition.

Let $A=(a_{ij})$, $i,j=1,\dots , n$, be a matrix over a ring. Let
$\sigma $ be a permutation of $\{1,\dots , n\}$. Write $\sigma $
as a product of disjoint cycles such that each cycle starts with a
smallest number. Since disjoint cycles commute, we may write
$$
\sigma =(k_{11}\dots k_{1j_1})(k_{21}\dots k_{2j_2})\dots (k_{m1}\dots k_{mj_m})
$$
where for each $i$, we have $k_{i1}<k_{ij}$ for all $j > 1$, and
$k_{11} > k_{21} \dots > k_{m1}$. This expression is unique. Let
$p(\sigma )$ be a parity of $\sigma $. Define Moore determinant
$M(A)$ by
$$
M(A)=\sum _{\sigma \in S_n}p(\sigma)a_{ k_{11}, k_{12}}\dots  a_{
k_{1j_1}, k_{11}} a_{ k_{21}, k_{22}}\dots a_{ k_{mj_m}, k_{m1}}.
\tag 4.1
$$

One can write cycles starting with the largest element. This will
lead to the same result.

If $A$ is Hermitian quaternionic matrix then $M(A)$ is a real number. Moore determinants
have nice features and are widely used (see, for example, {\cite {Al, As, Dy1}).

The following Theorem 4.3 shows that determinants of Hermitian
quaternionic matrices can be obtained through our general
approach. We prove first, that for quaternionic Hermitian matrix
$A$ determinants $D_{I,I}(A)$ coincide up to a sign.

Recall that $\Delta (A)=D_{I,I}(A)$ for $I=\{1,\dots , n\}$ and
that $\Delta (A)$ is a pre-Dieudonne determinant in the sense of
\cite {Dr}. If $A$ is Hermitian then $\Delta (A)$ is a product of
real numbers and, therefore, $\Delta (A)$ is real.

\proclaim{Proposition 4.2} Let $p(I)$ be the parity of the
ordering $I$. Then $\Delta (A)=p(I)p(J)D_{I,J}(A)$.
\endproclaim
The proof follows from homological relations for quasideterminants.

\proclaim {Theorem 4.3} Let $A$ be a Hermitian quaternionic
matrix. The $\Delta (A)=M(A)$.
\endproclaim
\demo {Proof} We use Sylvester formula for quasideterminants
(Theorem 1.7).

For $i,j=2,\dots , n$ set
$$
B_{ij}=\left ( \matrix
a_{11} & a_{1j}\\
a_{i1} & a_{ij}
\endmatrix \right ).
$$
$B_{ij}$ is a Hermitian matrix. Let $b_{ij}=M(B_{ij})$  and $c_{ij}=|B_{ij}|_{11}$.

Note that $B=(b_{ij})$ and $C=(c_{ij})$  also are Hermitian
matrices. It follows from (3.1) that  $M(A)=a_{nn}^{2-n}M(B)$.
Note, that $M(B)=a_{nn}^{n-1}M(C)$, therefore,  $M(A)=a_{nn}M(C)$.

By the noncommutative Sylvester formula (Theorem 1.7) ,
$|A|_{11}=|C|_{11}$, $|A^{11}|_{22}=|C^{11}|_{22}$, $\dots $. So,
$$
|A^{11}|_{22}|A^{11}|_{22}\dots |A^{12\dots n-1, 12\dots n-1}|_{n-1, n-1}=
|C^{11}|_{22}|C^{11}|_{22}\dots |C^{12\dots n-1, 12\dots n-1}|_{n-1, n-1}.
$$
The left-hand product equals  $\Delta (A)a_{nn}^{-1}$ and the
right-hand product equals
 $\Delta (C)$, so  $\Delta (A)= \Delta (C)a_{nn}=M(A)$.
\enddemo

\subhead 4.3. Moore determinants and norms of quaternionic
matrices \endsubhead

For a generic quaternionic matrix $A=(a_{ij})$, $i,j=1,\dots , n$,
define quaternionic norm $\nu (A)$ by induction:

i) $\nu (A)=\nu (a)$ if $A$ is the $1\times 1$-matrix $(a)$,

ii)  $\nu (A)=\nu (|A|_{11})\nu (A^{11})$.

It follows that $\nu (A)\in \Bbb R_{\geq 0}$.

\proclaim{Proposition 4.4} For generic matrices  $A, B$

i) $\nu (AB)=\nu (A)\nu (B)$,

ii) $\nu (A)=\Delta (A)\Delta (A^*)=\Delta (AA^*)$.
\endproclaim

Our proof is based on properties of quasideterminants.

Since $AA^*$ is a Hermitian matrix one has the following
\proclaim{Corollary 4.5} $\nu (A)=M(AA^*)$.
\endproclaim

It follows that $\nu (A)$ is a polynomial in the $a_{ij}$'s and
the $\bar a_{ij}$'s and so it can be defined for all matrices $A$
by continuity. The continuation posseses the multiplicative
property.

\subhead 4.4. Study determinants \endsubhead

In Section 1.1 a homomorphism $\theta :\Bbb H\rightarrow M(2, \Bbb
C)$ was defined. This homomorpism can be extended to homomorhism
of matrix algebras
$$\theta _n: M(n, \Bbb H) \rightarrow M(2n, \Bbb C).$$ Let
$A=(a_{ij})\in M(n, \Bbb H)$, set $\theta _n (A)=(\theta
(a_{ij}))$.

In 1920, E. Study \cite {S} defined a determinant $S(A)$ of  a
quaternionic matrix $A$ of order $n$ by setting $S(A)=\det \theta
_n (A)$. Here $\det $ is the standard determinant of a complex
matrix.

The following proposition is well known (see \cite {As}).
\proclaim{Proposition 4.6} For any quaternionic matrix $A$
$$S(A)=M(AA^*).$$
\endproclaim

The proof in \cite {As} was based on properties of eigenvalues of
quasideterminant matrices.
Our proof based on Sylvester identity and homological relations actually shows that $S(A)=\nu (A)$
for a generic matrix $A$.

\head V. Algebras with anti-involutions \endhead

All our results are valid for algebra $R$ with anti-involutions $I:R\rightarrow R$
such that subalgebra $R_0=\{r\in R: \ I(r)=r\}$ is central and the monoid of invertible elements
in $R$ is open in some topology.

Let $R$ be an algebra over a field $F$. An anti-involution of $R$
is a $F$-linear map $I: R\rightarrow R$ such that $I(xy)=I(y)I(x)$
and $I^2(x)=x$ for any $x,y \in R$. It is common to write $\bar
x$ instead of $I(x)$ for $x\in R$.

\enddocument

\bye